\documentclass[12pt,twoside]{article} 
\usepackage[english]{babel}
\usepackage[latin1]{inputenc} 
\usepackage{amsmath}
\usepackage{amssymb,amsfonts}
\usepackage{graphicx}                   
\usepackage{times,amssymb,amscd} 

\newcommand{\bA}{\mathbf{A}}

\newcommand{\bE}{\mathbf{E}}
\newcommand{\bG}{\mathbf{G}}
\newcommand{\bH}{\mathbf{H}}

\newcommand{\bL}{\mathbf{L}}

\newcommand{\bR}{\mathbf{R}}

\newcommand{\bS}{\mathbf{S}}
\newcommand{\bV}{\mathbf{V}}

\newcommand{\be}{\mathbf{e}}

\newcommand{\bx}{\mathbf{x}}
\newcommand{\by}{\mathbf{y}}

\newcommand{\bg}{\mathbf{g}}

\newcommand{\bI}{\mathbf{I}}

\newcommand{\BV}{\boldsymbol{V}}

\newcommand{\Be}{\boldsymbol{e}}
\newcommand{\Bu}{\boldsymbol{u}}
\newcommand{\Bv}{\boldsymbol{v}}

\newcommand{\cD}{\mathcal{D}}

\newcommand{\cP}{\mathcal{P}}
\newcommand{\cS}{\mathcal{S}}

\newcommand{\cB}{\mathcal{B}}

\newcommand{\EUC}{\bE^3}

\newcommand{\HYP}{\bH^3}
\newcommand{\SXR}{\bS^2\!\times\!\bR}
\newcommand{\HXR}{\bH^2\!\times\!\bR}
\newcommand{\SLR}{\widetilde{\bS\bL_2\bR}}
\newcommand{\NIL}{\mathbf{Nil}}
\newcommand{\SOL}{\mathbf{Sol}}

\begin{document}
\pagestyle{myheadings}
\markboth{\centerline{Jen\H o Szirmai}}
{A candidate to the densest packing with equal balls  $\dots$}
\title
{A candidate to the densest packing with equal balls in the Thurston geometries \footnote{AMS Classification 2010: 52C17, 52C22, 53A35, 51M20}}

\author{\normalsize Dedicated to Professor Emil {\sc Moln\'ar} \\  \normalsize on the Occasion of His 70th Birthday \vspace{4mm}\\
 Jen\H o Szirmai \\
% \thanks{This paper is supported by ....} \\
\normalsize Budapest University of Technology and \\ 
\normalsize Economics Institute of Mathematics, \\
\normalsize Department of Geometry \\
\normalsize szirmai@math.bme.hu\\
\date{\normalsize{\today}}}

%%%%%%%%%%%%%%%%%%%%%%%%%%%%%%%%%%%%%%%%%%%%
%\footnote{AMS Classification 2000: 52C17, 52C22, 53A35, 51M20}
%%%%%%%%%%%%%%%%%%%%%%%%%%%%%%%%%%%%%%%%%%%%

\maketitle
\begin{abstract}
 
The ball (or sphere) packing problem with equal balls, without any symmetry assumption, in a $3$-dimensional space of constant curvature 
was settled by B\"or\"oczky and 
Florian for the hyperbolic space $\HYP$ in \cite{BF64} and by proving the famous Kepler conjecture by Hales \cite{H} for the Euclidean space $\EUC$. 
The goal of this paper is to extend the problem of finding the densest geodesic 
ball (or sphere) packing for the other $3$-dimensional homogeneous geometries (Thurston geometries)
$$
\SXR,~\HXR,~\SLR,~\NIL,~\SOL, 
$$
where a transitive symmetry group of the  ball packing is assumed, one of the discrete isometry groups of the considered space.

Moreover, we describe a candidate of the densest geodesic ball packing. The greatest density until now is
$\approx 0.85327613$ that is not realized by packing with equal balls of the hyperbolic space $\HYP$. However,
it attains e.g. at horoball packing of $\overline{\bH}^3$ where the ideal centres of horoballs lie on the absolute figure of
$\overline{\bH}^3$  inducing the regular ideal simplex tiling $(3,3,6)$ by its Coxeter-Schl\"afli symbol.
In this work we present a geodesic ball packing in the $\SXR$ geometry whose density is $\approx 0.87499429$. The extremal configuration is described in
Theorem 2.8, Our conjecture and further remarks are summarized in Section 3.
\end{abstract}
%%%%%%%%%%%%%%%%%%%%%%%%%%%%%%%%%%%%%%%%%%%%

%%%%%%%%%%%%%%%%%%%%%%%%%%%%%%%%%%%%%%%%%%% 
\newtheorem{Theorem}{Theorem}[section]
\newtheorem{corollary}[Theorem]{Corollary}
\newtheorem{lemma}[Theorem]{Lemma}
\newtheorem{exmple}[Theorem]{Example}
\newtheorem{definition}[Theorem]{Definition}
\newtheorem{rmrk}[Theorem]{Remark}
\newtheorem{conjecture}[Theorem]{Conjecture}
\newtheorem{proposition}[Theorem]{Proposition}
\newenvironment{remark}{\begin{rmrk}\normalfont}{\end{rmrk}}
\newenvironment{example}{\begin{exmple}\normalfont}{\end{exmple}}
\newenvironment{acknowledgement}{Acknowledgement}

%%%%%%%%%%%%%%%%%%%%%%%%%%%%%%%%%%%%%%%%%%%%%%%%%%%%%%%%%%%%%%%%%%%% 

%============================================================================%
%                             the main article                               %
%============================================================================%

%%%%%%%%%%%%%%%%%%%%%%%%%%%%%%%%%%%%%%%%%%%%%%%%%%%%%%%%%%%%%%%%%%%%%%%%%%%%%%
\section{Preliminary results}
\subsection{Geodesic ball packings in spaces of constant curvature}
Finding the densest (not necessarily periodic) packing of balls in the 3-dimen\-sional Euclidean space is known as the Kepler Problem: 
{\it No packing of spheres of the same radius has a density greater than the face-centered cubic packing.}
This density can be realized by hexagonal layers (in continuum many ways). 
This conjecture was first published by Johannes Kepler in his monograph {\it The Six-Cornered Snowflake (1611)}, 
this treatise inspired by his correspondence with Thomas Harriot (see Cannonball Problem). 
In 1953, L\'aszl\'o Fejes T\'oth reduced the Kepler conjecture to an enormous calculation procedure that involved specific cases, 
and later suggested that computers might be helpful for solving the problem. In this way
the above four hundred year mathematical problem has finally been solved by  
Thomas Hales \cite{H}. He had proved that the guess of Kepler from 1611 was correct.
In mathematics sphere packing problems concern arrangements of non-overlapping equal spheres (rather balls) which fill a space. 
Space is the usual three-dimensional Euclidean space. However, ball (sphere) packing problems can be generalized to the other
$3$-dimensional Thurston geometries, but a difficult problem is -- similarly to the hyperbolic space -- the exact definition of the packing density.

In an $n$-dimensional space of constant curvature $\bE^n$, $\bH^n$, $\bS^n$ $(n\ge2)$ let $d_n(r)$ be the density of $n+1$ spheres
of radius $r$ mutually touching one another with respect to the simplex spanned by the centres of the spheres. L.~Fejes T\'oth and H.~S.~M.~Coxeter
conjectured that in an $n$-dimensional space of constant curvature the density of packing balls of radius $r$ can not exceed $d_n(r)$.
This conjecture has been proved by C.~Rogers in the Euclidean space $\bE^n$ \cite{Ro64}. 
The 2-dimensional spherical case was settled by L.~Fejes T\'oth in \cite{F65} and in \cite{B78} K.~B\"or\"oczky proved the following generalization:
\begin{Theorem}[K.~B\"or\"oczky]
In an $n$-dimensional space of constant curvature consider a packing of spheres of radius $r$. In spherical space suppose, that $r<\frac{\pi}{4}$.
Then the density of each sphere in its Dirichlet-Voronoi cell cannot exceed the density of $n+1$ spheres of radius $r$ mutually
touching one another with respect to the simplex spanned by their centres.
\end{Theorem}
\begin{rmrk}
\item[1.] In the hyperbolic space $\bH^3$ this result can be extended to $r=\infty$ \cite{B78} where the densest
horoball packing can be realized by different regular arrangements \cite{KSz}.
\item[2.] If we allow horoballs of different types at the various vertices of a totally asymptotic simplex and 
generalize the notion of the simplicial density function in $\mathbf{H}^n$, $(n \ge 2)$ then the B\"or\"oczky--Florian type 
density upper bound does not remain valid for the fully asymptotic simplices \cite{Sz11-3}, \cite{Sz12}. 
\end{rmrk}
\subsection{Geodesic ball packings under discrete isometry groups}
E. {Moln\'ar} has shown in \cite{M97}, that the homogeneous 3-spaces
have a unified interpretation in the projective 3-sphere $\mathcal{PS}^3(\bV^4,\BV_4, \mathbb{R})$. 
In our work we shall use this projective model of each Thurston geometry and
the Cartesian homogeneous coordinate simplex $E_0(\be_0)$,$E_1^{\infty}(\be_1)$,$E_2^{\infty}(\be_2)$,
$E_3^{\infty}(\be_3)$, $(\{\be_i\}\subset \bV^4$ \ \text{with the}  {\text{unit}} {\text{point}} $E(\be = \be_0 + \be_1 + \be_2 + \be_3 ))$ 
which is distinguished by an origin $E_0$ and by the ideal points of coordinate axes, respectively. 
Moreover, $\by=c\bx$ with $0<c\in \mathbb{R}$ (or $c\in\mathbb{R}\setminus\{0\})$
defines a point $(\bx)=(\by)$ of the projective 3-sphere $\cP \cS^3$ (or that of the projective space $\cP^3$ where opposite rays
$(\bx)$ and $(-\bx)$ are identified). 
The dual system $\{(\Be^i)\}\subset \BV_4$ describes the simplex planes, especially the plane at infinity 
$(\Be^0)=E_1^{\infty}E_2^{\infty}E_3^{\infty}$, and generally, $\Bv=\Bu\frac{1}{c}$ defines a plane $(\Bu)=(\Bv)$ of $\cP \cS^3$
(or that of $\cP^3$). Thus $0=\bx\Bu=\by\Bv$ defines the incidence of point $(\bx)=(\by)$ and plane
$(\Bu)=(\Bv)$, as $(\bx) \text{I} (\Bu)$ also denotes it. Thus the Thurston geometries can be visualized in the affine 3-space $\bA^3$
(so in $\bE^3$) as well.

Let $X$ be one of the remaining 5 Thurston geometries 
$$
\SXR,~\HXR,~\SLR,~\NIL,~\SOL
$$
where the geodesic curves are generally defined as having locally minimal arc length between their any two 
(near enough) points. 
The equation systems of the parametrized geodesic curves $\gamma(\tau)$ in our model can be determined by the 
general theory of Riemann geometry. {\it Then geodesic sphere and ball can be usually defined as follows below. 
We consider only those geodesic ball packings which are transitively generated by discrete groups of isometries of $X$ and
the density of the packing is related to its Dirichlet-Voronoi cells.}

\begin{definition}
The distance $d(P_1,P_2)$ between the points $P_1\in X$ and $P_2 \in X$ is defined by the arc length of the geodesic curve 
from $P_1$ to $P_2$.
\end{definition}
 \begin{definition}
 The geodesic sphere of radius $\rho$ (denoted by $S_{P_1}(\rho)$) with centre at the point $P_1$ is defined as the set of all points 
 $P_2$ in the space with the condition $d(P_1,P_2)=\rho$. Moreover, we require that the geodesic sphere is a simply connected 
 surface without selfintersection of space $X$.
 \end{definition}
 \begin{definition}
 The body of the geodesic sphere of centre $P_1$ and of radius $\rho$ in space $X$ is called geodesic ball, denoted by $B_{P_1}(\rho)$,
 i.e. $Q \in B_{P_1}(\rho)$ iff $0 \leq d(P_1,Q) \leq \rho$.
 \end{definition}
 In the following let $\Gamma$ be a fixed group of isometries of $X$. We 
 will denote by $d(P_1,P_2)$ the distance of two points $P_1$, $P_2$ (see Definition (1.3)). 
 \begin{definition}
 We say that the point set
 $$
 \cD(K)=\{P\in X\,:\,d(K,P)\leq d(K^\bg,P)\text{ for all }\bg\in\ \Gamma\}
 $$
 is the {Dirichlet--Voronoi cell} (D-V~cell) to $\Gamma$ around the kernel 
 point $K\in X$.
 \end{definition}
 \begin{definition}
 We say that
 $$
 \Gamma_P=\{\bg\in\Gamma\,:\,P^\bg=P\}
 $$
 is the \emph{stabilizer subgroup} of $P \in X$ in $\Gamma$.
 \end{definition}
 \begin{definition}
 Assume that the stabilizer $\Gamma_K=\bI$ the identity, i.e. $\Gamma$ acts simply transitively on 
 the $\Gamma$-orbit of $K \in X$. Then let $B_K$ denote the \emph{greatest ball} 
 of centre $K$ inside the D-V cell $\cD(K)$, moreover let $\rho(K)$ denote the 
 \emph{radius} of $B_K$. It is easy to see that
 $$
 \rho(K)=\min_{\bg\in\Gamma\setminus\bI}\frac12 d(K,K^\bg).
 $$
 \end{definition}
 \begin{definition}
 If the stabilizer $\Gamma_K > \bI$ then $\Gamma$ acts multiply transitively on 
 the $\Gamma$-orbit of $K \in X$. Then the greatest ball radius of $\cB_K$ is 
 $$
 \rho(K)=\min_{\bg\in\Gamma\setminus \Gamma_K}\frac12 d(K,K^\bg)
 $$
 where $K$ belongs to a 0- 1- or 2-dimensional region of $X$ (vertices, axes, reflection planes).
 \end{definition}
 In both cases the $\Gamma$-images of $B_K$ form a ball packing $\cB^\Gamma_K$ with centre 
 points $K^\bG$. 
 \begin{definition}
 The \emph{density} of ball packing $\cB^\Gamma_K$ is
 $$
 \delta(K)=\frac{Vol(B_K)}{Vol\cD(K)}.
 $$
 \end{definition}
 It is clear that the orbit $K^\Gamma$ and the ball packing $\cB^\Gamma_K$ have the 
 same symmetry group, moreover this group contains the starting 
 crystallographic group $\Gamma$:
 $$
 Sym K^\Gamma=Sym\cB^\Gamma_K\geq\Gamma.
 $$
 \begin{definition}\rm
 We say that the orbit $K^\Gamma$ and the ball packing $\cB^\Gamma_K$ is 
 \emph{characteristic} if $Sym K^\Gamma=\Gamma$, else the orbit is not 
 characteristic.
 \end{definition}
 \subsubsection{Simply transitive ball packings}
 
 Let $\Gamma$ be a fixed group of isometries in space $X$. 
 \emph{Our problem is} to find a 
 point $K\in\ X$ and the orbit $K^\Gamma$ for $\Gamma$ such that $\Gamma_K=\bI$ 
 and the density $\delta(K)$ of the corresponding ball packing 
 $\cB^\Gamma(K)$ is maximal. In this case the ball packing $\cB^\Gamma(K)$ is 
 said to be \emph{optimal.} 
 
 Our aim is to determine the maximal radius $\rho(K)$ of the balls, and the maximal density $\delta(K)$.
 The considered space groups could have free parameters, so we have to find the densest ball packing for fixed 
 parameters $p(\Gamma)$, then we have to vary them to get the optimal ball packing 
 \begin{equation}
 \delta(\Gamma)=\max_{K, \ p(\Gamma)}(\delta(K)). \tag{1.1}
 \end{equation}
 We look for the optimal kernel point
 in a 3-dimensional region, inside of a fundamental domain of $\Gamma$.
 
 \subsubsection{Multiply transitive ball packings}
 Similarly to the simply transitive case we have to find a kernel
 point $K\in\ X$ and the orbit $K^\Gamma$ for $\Gamma$ such that
 the density $\delta(K)$ of the corresponding ball packing 
 $\cB^\Gamma(K)$ is maximal but here $\Gamma_K \ne \bI$. This ball packing $\cB^\Gamma(K)$ is called \emph{optimal}, too. 
 In this multiply transitive case we look for the optimal kernel point $K$ 
 in possible 0- 1- or 2-dimensional regions $\mathcal{L}$, respectively. 
 Our aim is to deteremine the maximal radius $\rho(K)$ of the balls, and the maximal density $\delta(K)$.
 The considered space group can have also
 free parameters $p(\Gamma)$, then we have to find the densest ball packing for fixed 
 parameters, and vary them to get the optimal ball packing. 
 \begin{equation}
 \delta(\Gamma)=\max_{K\in \mathcal{L}, \ p(\Gamma)}(\delta(K)) \tag{1.2}
 \end{equation}
 \subsubsection{Packings in $\NIL$ space}
 
 {W. Heisenberg}'s famous real matrix group provides a non-commutative translation group of an affine 3-space.
 The $\NIL$ geometry, which is one of the eight homogeneous 
 Thurston 3-geometries, can be derived from this matrix group \cite{M97}. 
  
 In \cite{Sz07-2} I have investigated the geodesic balls of the $\NIL$ space and computed their volume, 
 introduced the notion of the $\NIL$ lattice, $\NIL$ parallelepiped and the density of the lattice-like ball packing.
 Moreover, I have determined the densest lattice-like geodesic ball packing by a type of $\NIL$ lattices. 
 The density of this packing is $\approx 0.78085$, may be surprising enough
 in comparison with the analogous Euclidean result $\frac{\pi}{\sqrt{18}} \approx 0.74048$. The kissing number of every ball
 in this packing is 14. {\it By my conjecture the densest geodesic ball packing belongs to the above ball arrangement in $\NIL$ space.}
 The symmetry group of this packing has also been described in \cite{MSz}.
 
 \subsubsection{Packings in $\HXR$ space}
 
 This Seifert fibre space is derived by the direct product of the hyperbolic plane $\bH^2$ and the real line $\bR$. 
 In \cite{Sz11} I have determined the geodesic balls of $\HXR$ space and computed their volume, 
 defined the notion of the geodesic ball packing and its density. 
 Moreover, I have developed a procedure to determine the density of the simply or multiply transitive geodesic ball packings for 
 generalized Coxeter space groups of $\HXR$ and apply this algorithm to them. 
 For the above space groups the Dirichlet-Voronoi cells are "prisms" in $\HXR$ sense. 
 The optimal packing density of the generalized Coxeter space groups is: $\approx 0.60726$. I'am sure, that in this space there are
 denser ball packings.
\begin{rmrk}
\begin{enumerate}
\item[1.] So far there are no results for the geodesic ball packings in $\SOL$ and $\SLR$ geometries (by my knowledge).
\item[2.] In $\NIL$ and $\SOL$ spaces I have studied the so-called {\it translation ball packings} (see \cite{Sz07-2}, \cite{Sz10})
but I do not consider these cases in this work. 
\end{enumerate}
\end{rmrk}
\section{On $\SXR$ space}
The $\SXR$ geometry can be derived by the direct product of the spherical plane $\bS^2$ and the real line $\bR$.
In \cite{F01} J.~Z. {Farkas} has classified and given the complete list of the space groups in $\SXR$.
The $\SXR$ manifolds up to similarity and diffeomorphism were classified by E. {Moln\'ar} and J.~Z. {Farkas} in \cite{FM01}.
In \cite{Sz11-1} I have studied the geodesic balls and their volumes in $\SXR$ space, moreover I have introduced the notion of geodesic ball packing 
and its density and determined the densest simply and multiply transitive geodesic ball packings for generalized Coxeter space groups of $\SXR$, respectively.
The density of the densest packing is $\approx 0.82445$.

In paper \cite{Sz11-2} I have studied the simply transitive locally optimal ball packings to the $\SXR$ space groups having Coxeter point groups 
and at least one of the generators is a non-trivial glide reflection. 
I have determined the densest simply transitive geodesic ball arrangements for the above space groups, moreover I computed their optimal 
densities and radii.
The density of the densest packing is $\approx 0.80408$.

Now, we shall discuss the simply and multiply transitive ball packings to a given space group. But let us start first with the necessary concepts.
The points in $\SXR$ geometry are described by $(P,p)$ where $P\in \bS^2$ and $p\in \bR$. 
The isometry group $Isom(\SXR)$ of $\SXR$ can be derived by the direct product 
of the isometry group of the spherical plane $Isom(\bS^2)$ and the isometry group of the real line $Isom(\bR)$. 

The structure of an isometry group $\Gamma \subset Isom(\SXR)$  is the following: $\Gamma:=\{(A_1 \times \rho_1), \dots (A_n \times \rho_n) \}$, where
$A_i \times \rho_i:=A_i \times (R_i,r_i):=(g_i,r_i)$, $(i \in \{ 1,2, \dots n \}$ and $A_i \in Isom(\bS^2)$, $R_i$ is either the identity map 
$\mathbf{1_R}$ of $\bR$ or the point reflection $\overline{\mathbf{1}}_{\mathbf{R}}$. $g_i:=A_i \times R_i$ is called the linear part of the transformation
$(A_i \times \rho_i)$ and $r_i$ is its translation part. 
The multiplication formula is the following:
\begin{equation}
(A_1 \times R_1,r_1) \circ (A_2 \times R_2,r_2)=((A_1A_2 \times R_1R_2,r_1R_2+r_2). \tag{2.1}
\end{equation}
\begin{definition}
A group of isometries $\Gamma \subset Isom(\SXR)$ is called {\it space group} if the linear parts form a finite group $\Gamma_0$ called the point group of 
$\Gamma$, moreover, the translation parts to the identity of this point group are required to form a one dimensional lattice $L_{\Gamma}$ of $\bR$.
\end{definition}
\begin{rmrk}
It can be proved, that the space group $\Gamma$ exactly described above has a compact fundamental domain $\mathcal{F}_\Gamma$.  
\end{rmrk}
We characterize the spherical plane groups by the {\it Macbeath-signature} (see \cite{M}, \cite{T}).

{\it In this paper we deal with a class of the $\SXR$ space groups {\bf 4q.~I.~2} (with a natural parameter $q \ge 2$, see \cite{F01}). 
Each of them belongs to the glide reflection groups, i.e. 
the generators $\bg_i, \ (i=1,2,\dots m)$ of its point group $\Gamma_0$ 
are reflections and at least one of the possible translation parts of the above generators differs from zero (see \cite{Sz11-2})}.

\subsection{Geodesic curves and balls in $\SXR$ space}

In \cite{Sz11-1} and \cite{Sz11-2} I have described the equation system of the geodesic curve and so the geodesic sphere: 
\begin{equation}
  \begin{gathered}
   x(\tau)=e^{\tau \sin{v}} \cos{(\tau \cos{v})}, \\ 
   y(\tau)=e^{\tau \sin{v}} \sin{(\tau \cos{v})} \cos{u}, \\
   z(\tau)=e^{\tau \sin{v}} \sin{(\tau \cos{v})} \sin{u},\\
   -\pi < u \le \pi,\ \ -\frac{\pi}{2}\le v \le \frac{\pi}{2} \tag{2.2}
  \end{gathered}
\end{equation}
of radius $\rho=\tau \ge 0$, centre $x(0)=1,~y(0)=0,~z(0)=0,$ and with longitude $u$, altitude $v$, as geographical coordinates.

In \cite{Sz11-1} I have proved that geodesic sphere $S(\rho)$ in $\SXR$ space is a simply connected surface
in $\mathbf{E}^3$ if and only if $\rho \in [0,\pi)$, because if $\rho\ge\pi$ then there is at least one $v \in 
[-\frac{\pi}{2},\frac{\pi}{2}]$ so that $y(\tau,v)=z(\tau,v)=0$, i.e. selfintersection would occur.
Thus we obtain the following
\begin{proposition}
The geodesic sphere and ball of radius $\rho$ exists in the $\SXR$ space if and only if $\rho \in [0,\pi).$
\end{proposition}
We have obtained (see \cite{Sz11-1}) the volume formula of the geodesic ball $B(\rho)$ of radius $\rho$ by 
the metric tensor $g_{ij}$ and by the Jacobian of (2.2) 
and a careful numerical Maple computation for given $\rho$ by the following integral:
\begin{Theorem}
\begin{equation}
\begin{gathered}
Vol(B(\rho))=\int_{V} \frac{1}{(x^2+y^2+z^2)^{3/2}}\mathrm{d}x ~ \mathrm{d}y ~ \mathrm{d}z = \\ = \int_{0}^{\rho} \int_{-\frac{\pi}{2}}^{\frac{\pi}{2}} 
\int_{-\pi}^{\pi} 
|\tau \cdot \sin(\cos(v)\tau)| ~ \mathrm{d} u \ \mathrm{d}v \ \mathrm{d}\tau = \\ =
2 \pi \int_{0}^{\rho} \int_{-\frac{\pi}{2}}^{\frac{\pi}{2}} |\tau \cdot \sin(\cos(v)\tau)| ~ \mathrm{d} v \ \mathrm{d}\tau. \tag{2.3}
\end{gathered}
\end{equation}
\end{Theorem}

The fundamental domain of the studied space groups can be combined as a fundamental domain of the spherical group with a part of a real line segment. 
This domain is called $\SXR$ {\it prism}. 
\begin{figure}[ht]
\centering
\includegraphics[width=6cm]{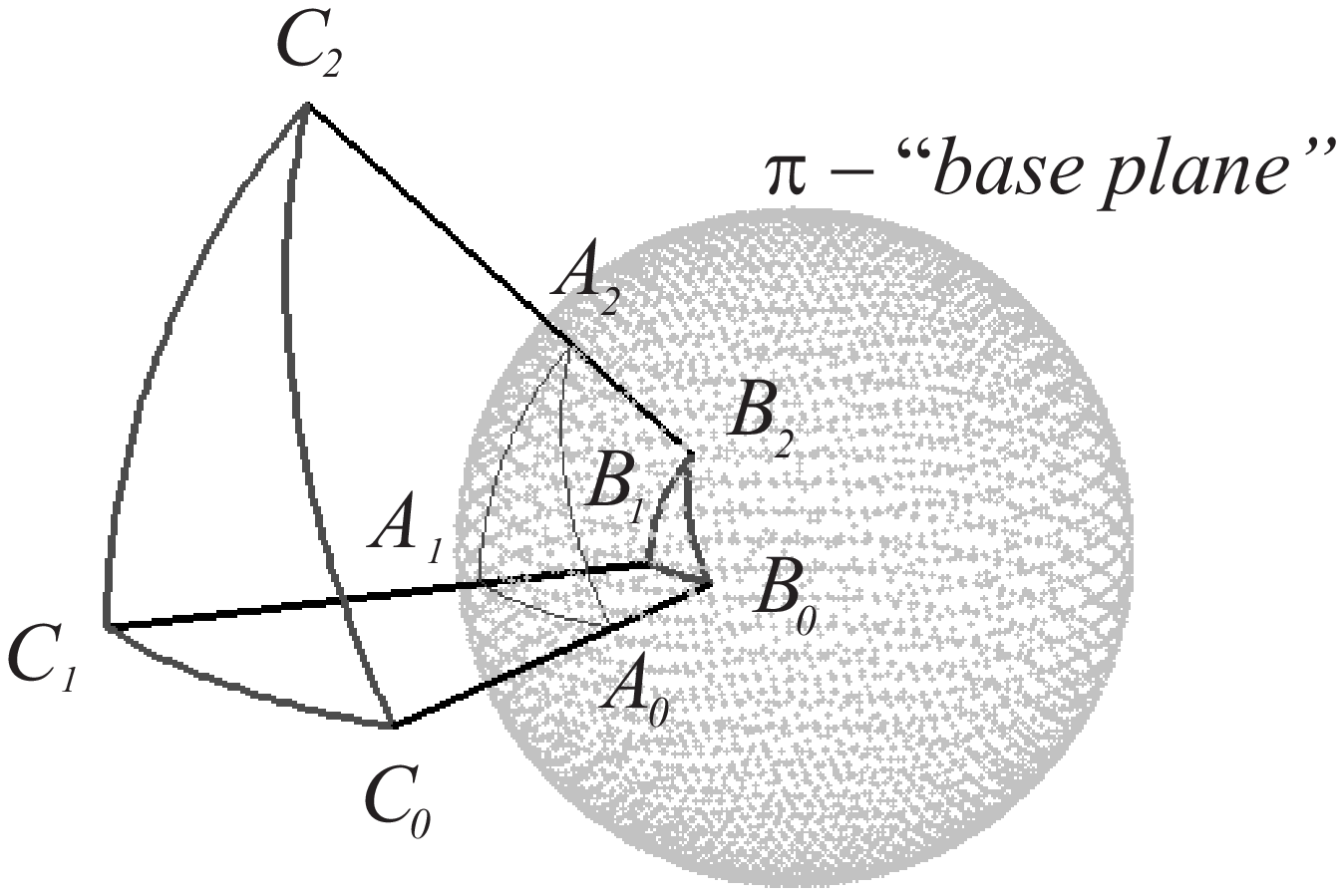} \includegraphics[width=4cm]{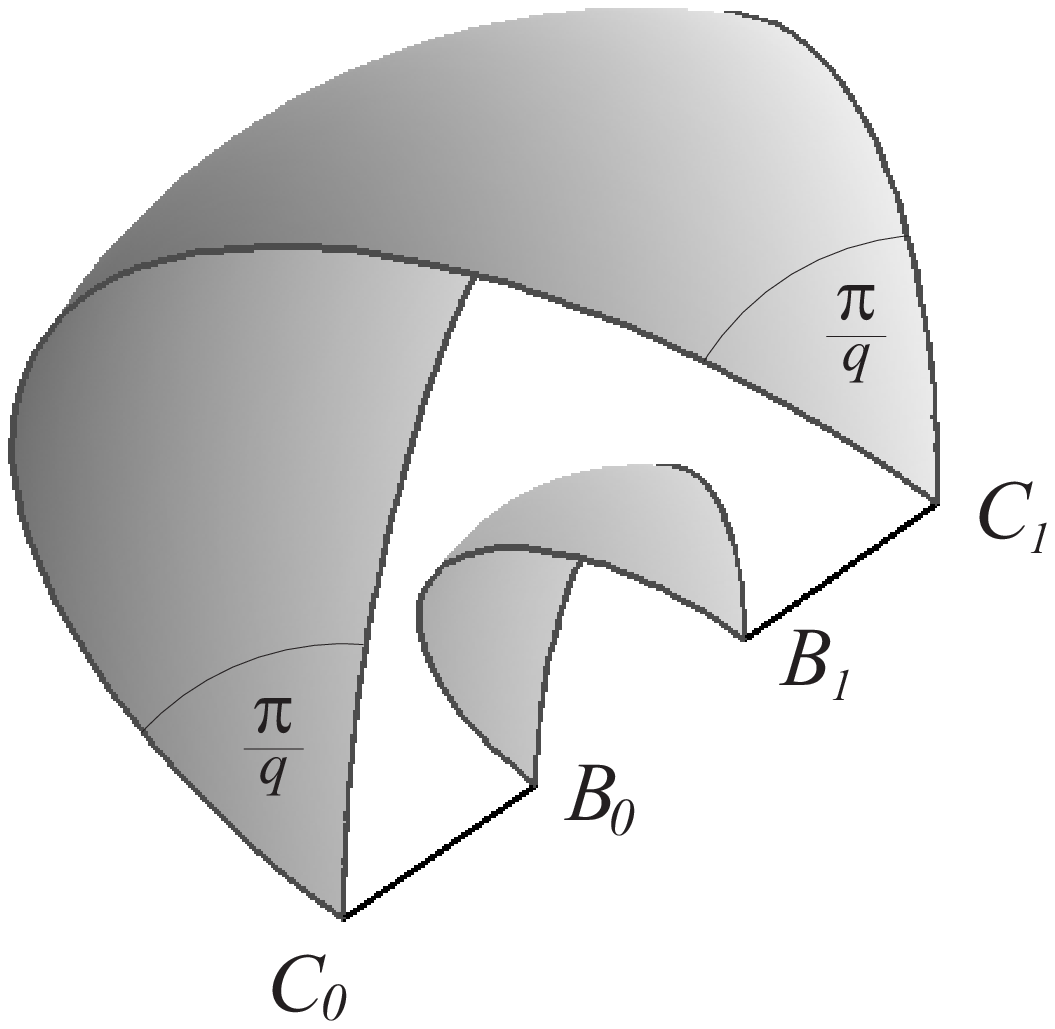}
\caption{Trigonal and digonal prisms}
\label{}
\end{figure}
In \cite{Sz11-1} we have shown the following 
\begin{Theorem}
The volume of a trigonal prism $\mathcal{P}_{B_0B_1B_2C_0C_1C_2}$ and of a digonal prism 
$\mathcal{P}_{B_0B_1C_0C_1}$ in $\SXR$ (see Fig.~1.a-b) can be computed by the following formula:
\begin{equation}
Vol(\mathcal{P})=\mathcal{A} \cdot h \tag{2.4}
\end{equation}
where $\mathcal{A}$ is the area of the spherical triangle $A_0A_1A_2$ or spherical digon $A_0A_1$ in the base plane $\Pi$ belonging to the fibre 
coordinate $t=0$ and
$h=B_0C_0$ is the height of the prism or digon, respectively.
\end{Theorem}
\subsection{Optimal ball packings for the space group {\bf 4q.~I.~2}, (q=2)}
We consider a $\SXR$ space group (see \cite{F01, Sz11-1, Sz11-2} 
with point group $\Gamma_0$ generated by three reflections $\bg_i, ~ (i=1,2,3)$ 
\begin{equation}
\begin{gathered}
(+,~0,~[~~]~ \{(2,2,q)\}),~ q \ge 2, \ \\ \Gamma_0=(\bg_1,\bg_2,\bg_3 - \bg_1^2,\bg_2^2,\bg_3^2,(\bg_1\bg_3)^2, (\bg_2\bg_3)^2), 
(\bg_1\bg_2)^q).
\end{gathered} \notag
\end{equation}
The possible translation parts $\tau_1$, $\tau_2$, $\tau_3$ of the corresponding generators of $\Gamma_0$ will be determined by (2.1) 
and by the defining relations of the above point group. 
Finally, we obtain six non-equivariant solutions
from the so-called Frobenius congruence relations:
$$(\tau_1,\tau_2,\tau_3) \cong (0,0,0),~\big(0,0,\frac{1}{2}\big),~\big(\frac{1}{2},\frac{1}{2},\frac{1}{2}\big),~
\big(\frac{1}{2},\frac{1}{2},0 \big),~\big(0, \frac{1}{2}, 0 \big),~\big(0, \frac{1}{2}, \frac{1}{2} \big).$$
If $(\tau_1,\tau_2,\tau_3) \cong ~(0,0,\frac{1}{2})$ then we have obtained the $\SXR$ space group {\bf 4q.~I.~2} (for a fixed $q$, $2\le q \in \mathbf{N}$).

The fundamental domain of the point group of the considered space group is a spherical triangle $A_1A_2A_3$ with 
angle $\frac{\pi}{q}$, $\frac{\pi}{2}$, $\frac{\pi}{2}$ in the base plane $\Pi$.
It can be assumed, that the fibre coordinate of the center of the optimal ball is zero and it is an 
interior point of $A_1A_2A_3$ triangle (see Fig.~2).  

{\it In the following we consider ball packings, only to $q=2$.}
\begin{figure}
\centering
\includegraphics[width=10cm]{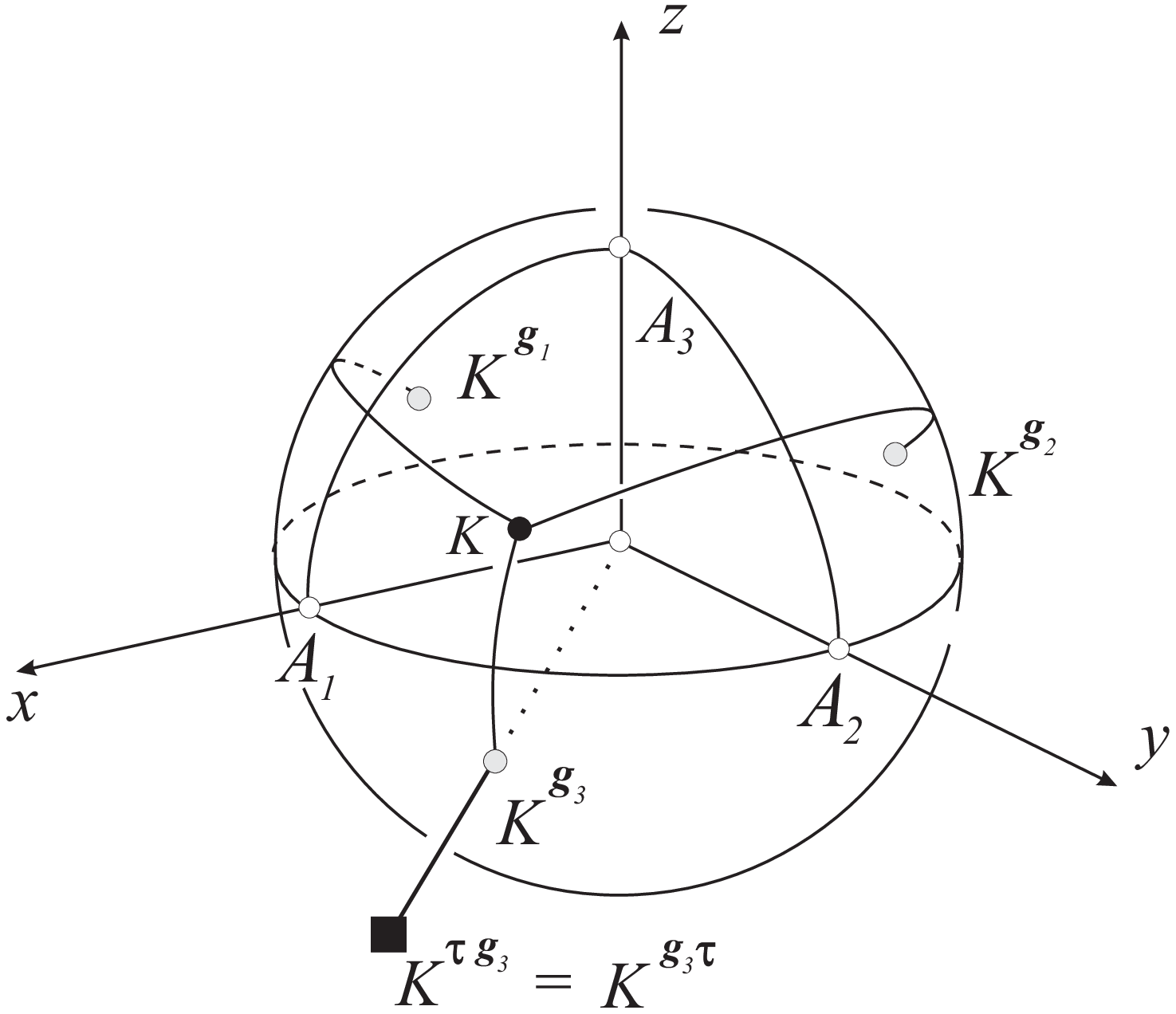}
\caption{}
\label{}
\end{figure}
We shall apply the above introduced Cartesian homogeneous coordinate system (see Fig.~2) and the usual geographical coordinates 
$(\phi, \theta), ~ (-\pi < \phi \le \pi, ~ -\frac{\pi}{2}\le \theta \le \frac{\pi}{2})$ 
of the sphere with the fibre coordinate $t \in \bR$.
\subsubsection{The simply transitive ball packing}

We consider an arbitrary interior point $K(x^0,x^1,x^2,x^3)=K(\phi,\theta)$ of spherical triangle $A_1A_2A_3$  
in the above coordinate system in our model (see Fig.~2).
Our aim is to determine the maximal radius $R$ of the balls, and the maximal density $\delta(K,R)$ where $K$ is an inner point of $A_1A_2A_3$ triangle
or $K\in A_1A_2$.

Thus, $q=2$, $\tau_3=\tau$, the $\SXR$ group is ${\bf 42.I.2}$ now, the $\bR$-translation is $2\tau$. 
(Instead of the above 6 cases we have only 4 space groups for the corresponding point group.)
The fundamental domain of the point group of the considered space group is a spherical triangle $A_1A_2A_3$ with 
angles $\frac{\pi}{2}$, $\frac{\pi}{2}$, $\frac{\pi}{2}$ in the base plane $\Pi$ (see Fig.~2).
It can be assumed by the homogeneity of $\SXR$, that the fibre coordinate of the center of the optimal ball is zero.
\begin{equation}
x^0=1, \ \ x^1= \cos{\phi} \cos{\theta},  \ \ x^2=\sin{\phi} \cos{\theta},  \ \ x^3= \sin{\theta} \tag{2.5}
\end{equation}

Let $\mathcal{B}_\Gamma(R)$ denote a geodesic ball packing of $\SXR$ with balls $B_K(R)$ of radius $R$ (to be determined yet) where their 
centres $K$ give rise to the orbit $K^\Gamma$. In the following we consider to each ball packing the possible smallest 
translation part $\tau(K,R)$ (see Fig.~2) depending on $\Gamma$, $K$ and $R$. 
A fundamental domain of $\Gamma$ $\widetilde{\mathcal{F}}_{\tau(K,R)}$ is the corresponding Dirichlet-Voronoi cell where its volume 
is equal to the volume of the prism which is given by the fundamental domain of the point group $\Gamma_0$ of $\Gamma$ and by the
height $2\tau(K,R)$.
(It is clear that the optimal ball ${B}_K$ has to touch some faces of its D-V cell to $\Gamma$ around the kernel 
point $K$.)
The images of $\widetilde{\mathcal{F}}_{\tau(K,R)}$ 
by discrete isometry group $\Gamma$ covers the $\SXR$ space without overlap. 
For the density of the packing it is sufficient to relate the volume of the optimal ball
to that of the solid $\widetilde{\mathcal{F}}_{\tau(K,R)}$.  
Analogously to the Euclidean case the density $\delta(R,K)$ of the geodesic
ball packing $\mathcal{B}_\Gamma(R)$ can be defined :
\begin{definition}
\begin{equation}
\delta(R,K):=\frac{Vol(\mathcal{B}_\Gamma(R) \cap \widetilde{\mathcal{F}}_{\tau(K,R)})}{Vol(\widetilde{\mathcal{F}}_{\tau(K,R)})}. \tag{2.6}
\end{equation}
\end{definition}
\begin{figure}[ht]
\centering
\includegraphics[width=6cm]{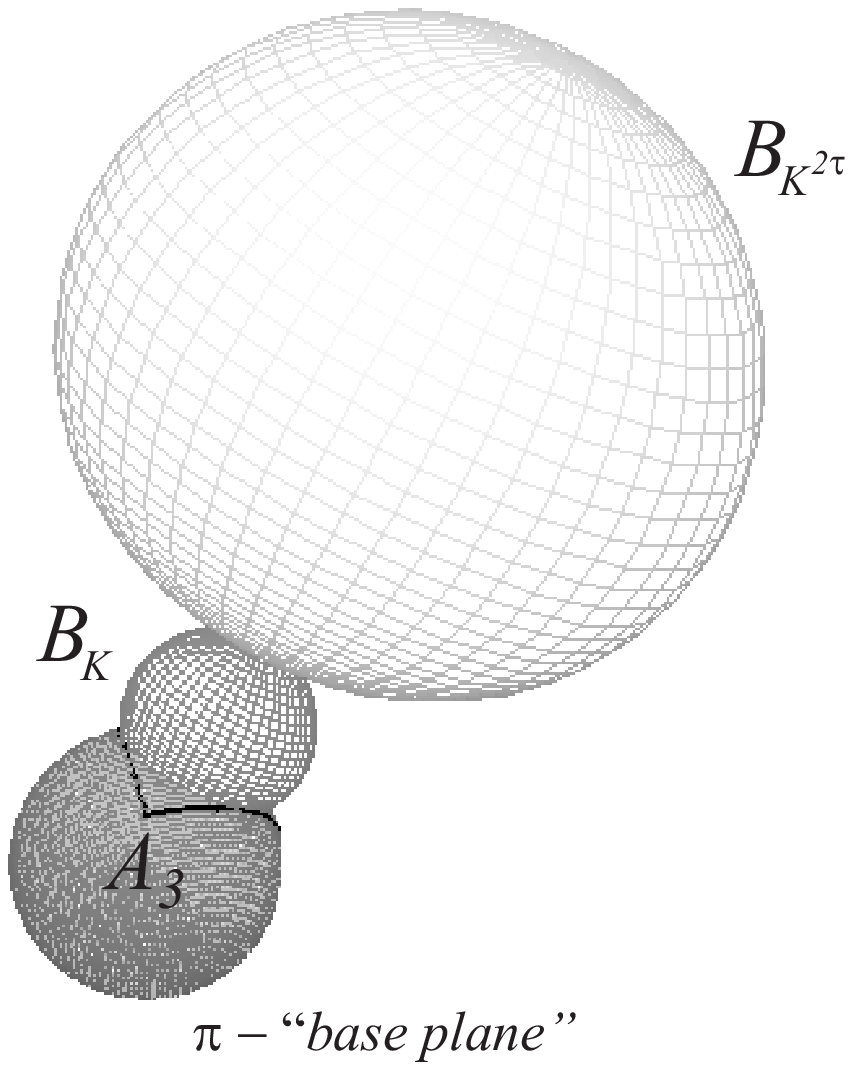} \includegraphics[width=6cm]{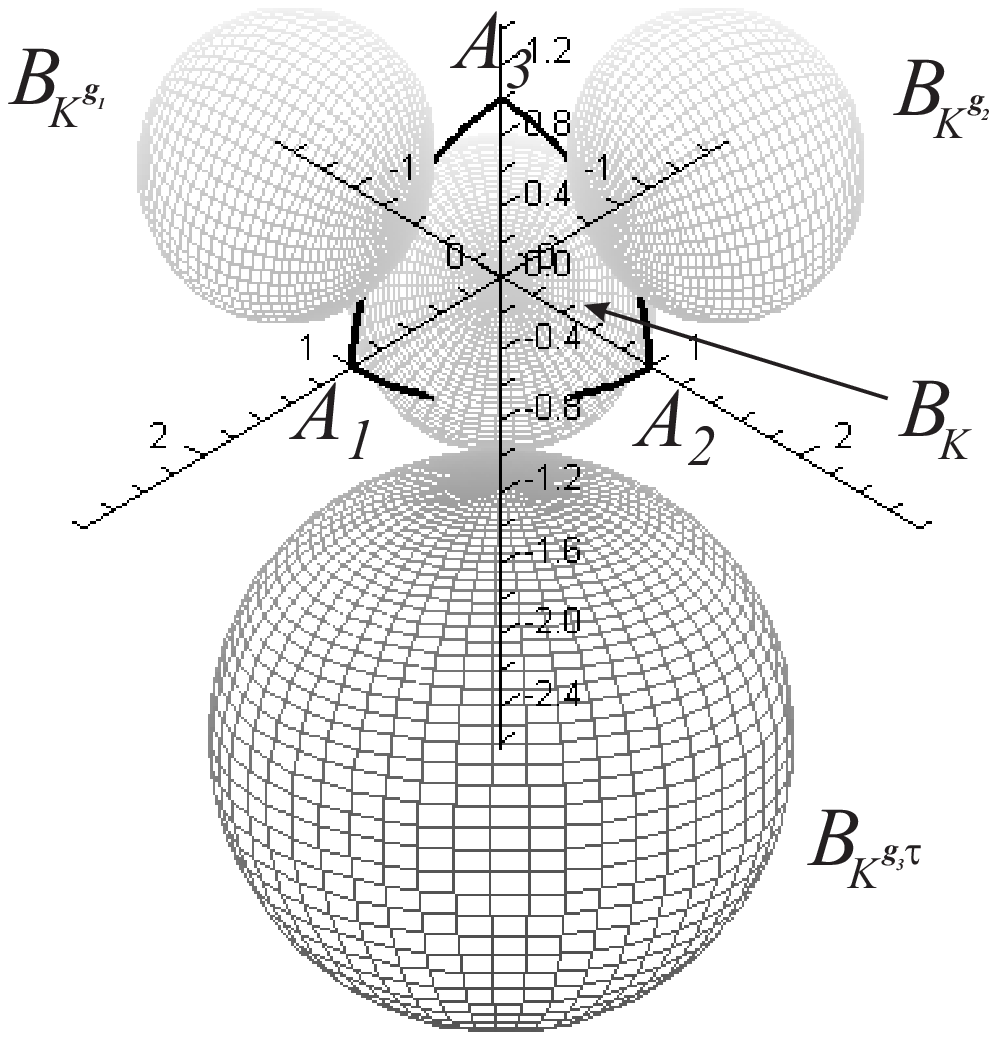}

a. \hspace{6cm} b.
\caption{Balls of the simply transitive optimal ball arrangement with base plane $\pi$ which is represented by the spherical triangle $A_1A_2A_3$ 
in figure $b$.}
\label{}
\end{figure}
\begin{figure}[ht]
\centering
\includegraphics[width=10cm]{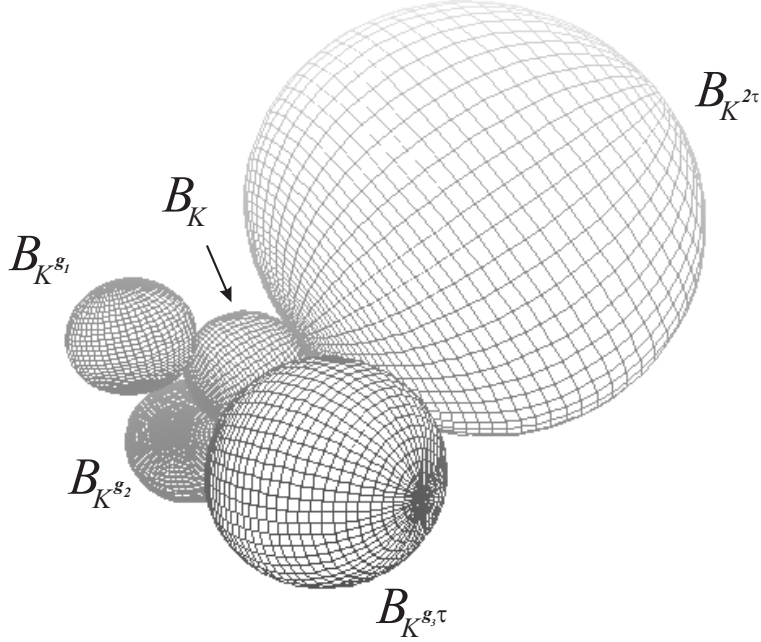} 
\caption{Some neighbouring balls of the simply transitive optimal ball arrangement $\mathcal{B}_{opt}(K,R)$=$\mathcal{B}_{\Gamma}(R_1,K_1)$}
\label{}
\end{figure}
The optimal ball arrangement $\mathcal{B}_{opt}(K,R)$ has to satisfy the following equations:
\begin{equation}
\begin{gathered}
(a) \ \ d(K,K^{g_1})=2R=d(K,K^{g_2}), \\
(b) \ \ d(K,K^{\tau g_3})=2R \le d(K,K^{\tau g_3 \tau g_3}) \le 4R=d(K,K^{2\tau}). \\
\end{gathered} \tag{2.7}
\end{equation}
We consider two main ball arrengements: 
\begin{enumerate}
\item We denote by $\mathcal{B}_{\Gamma}(R_1,K_1)$ those packing where requirements (2.7)  
and \\ $d(K,K^{2\tau}) = 2R$ hold (see Fig.~2).
\item We denote by $\mathcal{B}_{\Gamma}(R_2,K_2)$ those packing where requirements (2.7)  
and  \\ $d(K,K^{\tau g_3 \tau g_3})=d(K,K^{2\tau}) = 4R$ hold (see Fig.~2).
\end{enumerate}
First we determine the coordinates of the points $K_i$, ($i=1,~2$) ($K_i$ with parameters $\phi_i$ and $\theta_i$), the radius $R_i$ of the ball, 
the volume of a ball $B(R_i)$  and the density of the packing in both main cases. 
We get the following solutions by systematic approximation, where the computations were carried out by {\it Maple V Release 10} up to 30 decimals:
\begin{equation}
\begin{gathered}
\phi_1 = \frac{\pi}{4} \approx 0.78539816,\ \ \theta_1 \approx 0.55737781, \ \ R_1 \approx 0.64360446, \\
Vol(B(R_1)\approx 1.08624788, \ \ \delta(R_1,K_1) \approx 0.53722971.
\end{gathered} \tag{2.8}
\end{equation}
\begin{equation}
\begin{gathered}
\phi_2 = \frac{\pi}{4}\approx 0.78539816,\ \ \theta_2=0, \ \ R_2 = \frac{\pi}{4}\approx 0.78539816, \\
Vol(B(R_2)\approx 1.94735865, \ \ \delta(R_2,K_2) \approx 0.39461737.
\end{gathered} \tag{2.9}
\end{equation}
We obtain by careful discussion of the density function 
$\delta(R,K)$ ($R \in [R_1,R_2]$) of the considered ball packing  the following:
\begin{Theorem}
The ball arrangement $\mathcal{B}_{\Gamma}(R_1,K_1)$ (see Fig.~3.a-b and Fig.~4) provides the densest symply transitive ball packing belonging to the $\SXR$ 
space group {\bf 4q.~I.~2} $(q=2)$. 
\end{Theorem}

\subsubsection{The multiply transitive ball packings}
To determine the optimal multiply ball packing we have to study 2-cases:
\begin{enumerate}
\item[1.] $K$ is an inner point of the spherical geodesic segment $A_2A_3$ (or $A_1A_3$).
In this situation the point $K$ and its images by $\Gamma={\bf 4q.~I.~2}$ $(q=2)$, as the centers of the optimal ball arrangement 
$\mathcal{B}_{opt}(K,R)$ have to hold the following requirements because of an arbitrary ball of the optimal packing is fixed by 
its neighbouring balls:
\begin{equation}
\begin{gathered}
(a) \ \ d(K,K^{g_1})=2R=d(K,K^{\tau g_3}), \\
(b) \ \ 2R \le d(K,K^{\tau g_3 \tau g_3}) = d(K,K^{2 \tau}) \le  4R. \\
\end{gathered} \tag{2.10}
\end{equation}
It is easy to see that in this case we get the optimal packing if $K=A_2$ (or $K=A_1$) with the following data:
\begin{equation}
\begin{gathered}
\phi_3= \frac{\pi}{2} \approx 1.57079633,\ \ \theta_3 =0, \ \ R_3 = \frac{\pi}{2}\approx 1.57079633, \\
Vol(B(R_3))\approx 13.74539472, \ \ \delta(R_3,K_3) \approx 0.69634983.
\end{gathered} \tag{2.11}
\end{equation}
\item[2.] $K=A_3$.
The Fig.~5 shows the orbit of the point $K=A_3$ by the considered space group. The images of $K$ lie on a line through the origin and $A_3$.   
\begin{figure}[ht]
\centering
\includegraphics[width=9cm]{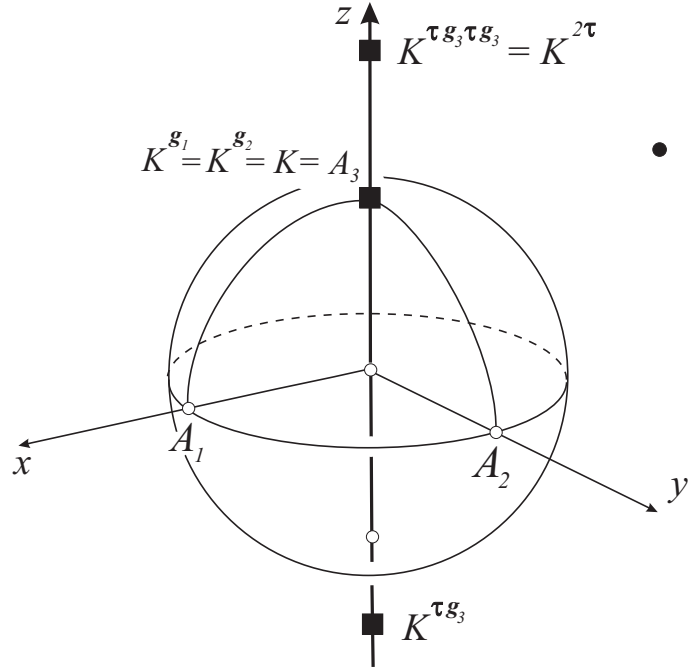}
\caption{The orbit of $K=A_3$ by the group $\Gamma={\bf 4q.~I.~2}$ $(q=2)$.}
\label{}
\end{figure}
Fig.~6 and Fig.~7. show this interesting ball arrangement whose data are the following

\begin{equation}
\begin{gathered}
\phi_4= \frac{\pi}{4} \approx 0.78539816,\ \ \theta_4 =\frac{\pi}{2} \approx 1.57079633, \ \ R_4 \approx 1.81379936, \\
Vol(B(R_4))\approx 20.00238509, \ \ \delta(R_4,K_4) \approx 0.87757183.
\end{gathered} \tag{2.12}
\end{equation}
The "outwards transformed" images of a balls sourround the previous balls (see Fig.~7.) thus the touching number of this packing is 4.   
Finally, we get the following
\begin{Theorem}
The ball arrangement $\mathcal{B}_{opt}(R_4,K_4)$ provides the densest multiply transitive ball packing belonging to 
the $\SXR$ space group {\bf 4q.~I.~2} $(q=2)$.
\end{Theorem}
\end{enumerate}
\begin{figure}[ht]
\centering
\includegraphics[width=9cm]{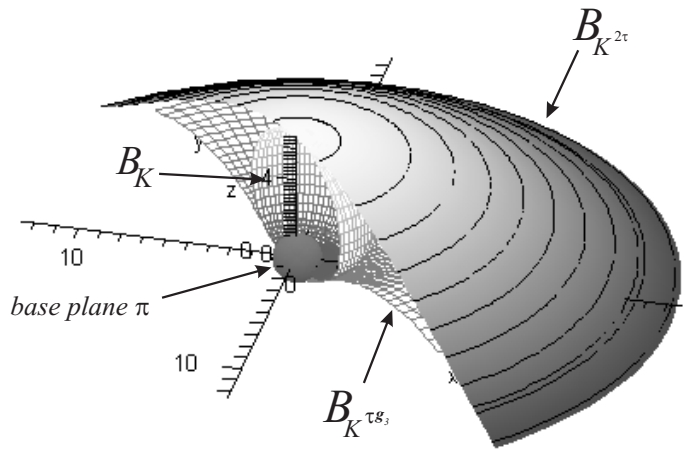}
\caption{The densest ball packing $\mathcal{B}_{opt}(R_4,K_4)$ is illustrated with parts of some neighbouring half spheres 
$B_K$, $B_{K^{\tau \bg_3}}$, $B_{K^{2\tau}}$  and the "base plane". The spheres $B_{K^{\tau \bg_3}}$, $B_{K^{2\tau}}$ sourround and touch the 
sphere $B_K$. The centres of the above spheres are described in Fig.~5.}  
\label{}
\end{figure}
\begin{figure}[ht]
\centering
\includegraphics[width=8cm]{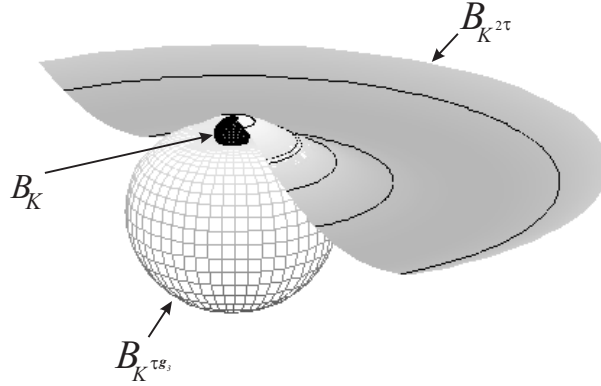} 
\caption{The densest ball packing is described by its balls $B_K$, $B_{K^{\tau g_3}}$ and a part of the shere $B_{K^{2\tau}}$.}
\label{}
\end{figure}
\section{The conjecture for the densest ball arrangement in Thurston geometries}
The notion of the density of an arbitrary congruent geodesic ball packing in spaces of constant curvature is known (see Section 1).

In Section 2 we have introduced the density function of the geodesic ball packings generated by a discrete isometry group in a given Thurston geometry.
This density is related to the Dirichlet-Voronoi cells generated by centres of the balls. 
For these ball packings we can formulate the following
\begin{conjecture}
Let $\cB$ be an arbitrary congruent geodesic ball packing in a Thurston geometry $X$, where $\cB$ generated by a discrete isometry group of $X$.
The above determined ball arrangement $\mathcal{B}_{opt}(R_4,K_4)$ with density $\delta(R_4,K_4) \approx 0.87757183$ provides 
the densest congruent geodesic ball packing in the Thurston geometries. 
\end{conjecture}
\begin{rmrk}
So far the shape of Dirichlet-Voronoi cells moreover
the equidistant surfaces of two points are unknown in some spaces (see \cite{PSSz-1}, \cite{PSSz-2}, \cite{PSSz-3}).
\end{rmrk}
The general definition of the density of congruent geodesic ball packings in the Thurston geometries is not settled yet, but by our investigation 
for any "good" density definition can be formulated the next   
\begin{conjecture}
The densest congruent geodesic ball packing in the Thurston geometries is realized by
the above ball arrangement $\mathcal{B}_{opt}(R_4,K_4)$ with density $\delta(R_4,K_4) \approx 0.87757183$. 
\end{conjecture}
We are working in similar problems in $\SLR$ and $\SOL$ space, too. 

In this paper we have mentioned only some problems in discrete geometry of Thurston spaces, but we hope that from these 
it can be seen that our projective
method suits to study and solve similar problems (\cite{Sz10}, \cite{Sz07-2}, \cite{Sz11-4}).

\end{document}